\documentclass[a4paper,12pt]{article}
\usepackage{amsmath,amsthm,amsfonts,latexsym,amscd}

\scrollmode

\swapnumbers
\theoremstyle{plain}
\newtheorem{theorem}{Theorem}[section]

\newtheorem{corollary}[theorem]{Corollary}
\newtheorem{proposition}[theorem]{Proposition}

\theoremstyle{definition}
\newtheorem{definition}[theorem]{Definition}
\newtheorem{example}[theorem]{Example}

\theoremstyle{remark}

\newcommand{\R}{\mathbb{R}}
\newcommand{\C}{\mathbb{C}}

\newcommand{\Z}{\mathbb{Z}}

\newcommand{\HH}{\mathbb{H}}

\newcommand{\abs}[1]{\lvert#1\rvert} 

\newcommand{\into}{\hookrightarrow}

\newcommand{\iso}{\cong}

\DeclareMathOperator{\ind}{ind}

{\catcode`@=11\global\let\c@equation=\c@theorem}

\begin{document}

\title{A counterexample to the (unstable) Gromov-Lawson-Rosenberg conjecture}
\author{Thomas Schick\thanks{
e-mail: thomas.schick@uni-math.gwdg.de\hspace{1cm}
www:~http://www.uni-math.gwdg.de/schick/
}\\Fachbereich Mathematik ---
Universit{\"a}t G\"ottingen\\
Bunsenstr.~3\\
 37073 G\"ottingen\\ Germany 
}
         
\maketitle

\begin{abstract}
Doing surgery on the $5$-torus, we construct a $5$-dimensional closed 
spin-manifold $M$ with $\pi_1(M)\iso \Z^4\times \Z/3$, so that the index
invariant in the $KO$-theory of the reduced $C^*$-algebra of $\pi_1(M)$
 is zero. Then we use the theory of minimal surfaces of Schoen/Yau to
show that this manifolds cannot carry a metric of positive scalar curvature.
The existence of such a metric is predicted by the (unstable)
 Gromov-Lawson-Rosenberg
conjecture.
\end{abstract}

\section{Obstructions to positive scalar curvature}
\begin{definition}
  A manifold $M$ which admits a metric of positive scalar curvature is called
a pscm-manifold.
\end{definition}

We start with a discussion of
 the index obstruction for spin manifolds to be pscm, 
constructed by
Lichnerowicz \cite{Lichnerowicz(1963)}, Hitchin \cite{Hitchin(1974b)} and in
the following refined version due to Rosenberg \cite{Rosenberg(1987)}.
\begin{theorem}
  Let $M^m$ be a closed spin-manifold, $\pi:=\pi_1(M)$. One can construct a
homomorphism, called index, from the singular spin bordism 
$\Omega_*^{spin}(B\pi)$ to the (real) $KO$-theory of the reduced real
$C^*$-algebra of $\pi$:
\[  \ind: \Omega_*^{spin}(B\pi)\to KO_*(C^*_{red}\pi) \]
Let $u:M\to B\pi$ be the classifying map for the universal covering of $M$.
If $M$ is pscm, then
\[ \ind([u:M\to B\pi]) =0\in KO_m(C^*_{red}) \]
\end{theorem}
Gromov/Lawson \cite{Gromov-Lawson(1983)} and Rosenberg \cite{Rosenberg(1983)}
 conjectured that the vanishing of the index should
also be sufficient for existence of a metric with positive scalar curvature
 on $M$ if $m\ge 5$.
 This was proven by Stefan Stolz \cite{Stolz(1992)}
 for $\pi=1$, and subsequently by him and other authors
 also for a few other groups 
\cite{Rosenberg(1983),Kwasik-Schultz(1990),%
Botvinnik-Gilkey-Stolz,Rosenberg-Stolz(1995)}.

In dimension $\ge 5$ there is only one additional obstruction for pscm,
the minimal surface method of  Schoen and Yau, which we will recall now.
(In dimension $4$, the Seiberg-Witten theory yields additional obstructions).
The first 
theorem is the differential geometrical backbone for the application of
minimal surfaces to the pscm problem:
\begin{theorem}\label{diffgeo}
Let $(M^m, g)$ be a manifold with positive scalar curvature, $\dim M=m\ge 3$.
 If $V$ is a smooth $m-1$-dimensional
submanifold of $M$ with trivial normal bundle,
 and if $V$ is a local minimum of the
volume functional, then $V$ admits a metric of positive scalar curvature.
\end{theorem}
\begin{proof}
Schoen/Yau: \cite[5.1]{Schoen-Yau(1979a)} for $m=3$, 
\cite[proof of Theorem1]{Schoen-Yau(1979)} for $m>3$.
\end{proof}

The next statement from geometric measure theory implies applicability of the
previous theorem if $\dim(M)\le 7$.
\begin{theorem}\cite[chapter 8]{Morgan(1988)} and references therein,
in particular \cite[5.4.18]{Federer(1969)}\label{regul}\\
  Suppose $M^m$ is a smooth orientable closed manifold, $\dim M=m\le 7$.
 Suppose
$0\ne x\in H_{m-1}(M,\Z)$. Then a smooth orientable
closed $m-1$-dimensional submanifold 
$V$
of $M$ exists which represents $x$ and which has minimal volume under all
currents which represent $x$. In particular, $V$ is a local minimum of the
volume functional with orientable (hence trivial) normal bundle.
\end{theorem}

This implies the following statement about homology and cohomology which was
observed by Stephan Stolz.
Let $X$ be any space.
\begin{definition}
  For $m\ge 2$ we define
\[ H_m^+(X,\Z):=\{ f_*[M]\in H_m(X,\Z);\; f:M^m\to X
\text{ and $M$ is pscm}\} \]
\end{definition}

\begin{corollary}\label{Hom_cond}
Let $X$ be any space, $\alpha\in H^1(X,\Z)$. Cap-product with $\alpha$
induces a map 
\[ \alpha\cap:H_m(X,\Z)\to H_{m-1}(X,\Z).\]
 If $3\le m\le 7$, then
$\alpha\cap$ maps $H_m^+(X,\Z)$ to $H_{m-1}^+(X,\Z)$.
\end{corollary}
\begin{proof}
If $f:M^m\to X$ represents $x\in H_m^+(X,\Z)$ and $M$ is a pscm-manifold,
 then by theorems \ref{diffgeo}
and \ref{regul}
the class $f^*\alpha\cap [M]$ is represented by $N^{m-1}\stackrel{j}{
\into} M$ and $N$ is also pscm.
 In particular $\alpha\cap x=f_*(f^*\alpha\cap[M] )$
is represented by $f\circ j:N\to X$, i.e.\ $\alpha\cap x\in H_{m-1}^+(X,\Z)$.
\end{proof}

\section{The Counterexample}
To produce a counterexample to the unstable Gromov-Lawson-Rosenberg
conjecture, we use the only  other known obstruction for
pscm, namely the minimal surface method explained above.

The fundamental group will be $\pi:= \Z^4\times\Z/3$. We start with the 
computation of
the $KO$-theory of the $C^*\pi$-algebra.
Note that the reduced $C^*$-algebra of the product of two groups is the 
(minimal)
tensor product of the individual $C^*$-algebras 
\cite[p.\ 14--15]{Schroeder(1993)}. By 
\cite[p.\ 14 and 1.5.4]{Schroeder(1993)}
\[ KO_n(C^*_{red}(\Z^4\times \Z/3))\cong 
\bigoplus_{i\in I} KO_{n-n_i}(C^*_{red}
(\Z/3));\qquad \abs{I}<\infty \]
For finite groups, it is well known that their $KO$-theory is a direct sum
of copies of the $KO$-theories of $\R$, $\C$ and $\HH$. In particular, it is
a direct sum of copies of $\Z$ and $\Z/2$. Therefore, the same is true for
$\pi$:
\begin{proposition}
  $KO_*(C^*_{red}\pi)$ is a direct sum of copies of $\Z$ and $\Z/2$. In
particular, its torsion is only $2$-torsion.
\end{proposition}

We will now construct a spin manifold $M^5$ with $\pi_1(M)=\pi$,
 so that the class
$[u:M\to B\pi]\in \Omega_5^{spin}$ is $3$-torsion. Then, automatically
\[ \ind(u:M\to B\pi)=0\in KO(C^*_{red}\pi) \]

\begin{example}
Let $p:S^1\to B\Z/3$ be a map so that $\pi_1(p)$ is surjective and equip $S^1$
with the spin structure induced from $D^2$.
This is $3$-torsion since
 $\tilde\Omega_1^{spin}(B\Z/3)\cong H_1(B\Z/3,\Z)\cong\Z/3$ (use
the Atiyah-Hirzebruch spectral sequence). Consider the
singular manifold 
\[ f = id\times p: S^1\overbrace{\times\cdot\times}^4 S^1\times S^1 \to
 S^1\overbrace{\times\cdot\times}^4 S^1\times B\Z/3 = B\pi \]
This is then $3$-torsion in $\Omega_5^{spin}(B\pi)$. Doing surgery
we can construct a bordism $F:W\to B\pi$ 
 in $\Omega_5^{spin}(B\pi)$ from $f$ to some $u: M\to B\pi$
where $u$ is an isomorphism on $\pi_1$. 
\end{example}
Now, $M$ is a manifold with trivial
index, and we have to show that it is not pcsm. Assume that the converse is
true.

We study the homology and cohomology of $\pi$ first.
By the K{\"u}nneth theorem
\[\begin{split}
  H_1(B\pi,\Z) &= x_1\Z\oplus \dots \oplus x_4\Z \oplus y\Z/3\\
  H^1(B\pi,\Z) &= a_1\Z\oplus \dots \oplus a_4\Z\\
  0\ne w= x_1\times \dots\times x_4\times y &\in H_5(B\pi,\Z)\\
  0\ne z= x_4\times y =a_1\cap (a_2\cap (a_3 \cap w))&\in H_2(B\pi,\Z)
\end{split} \]

We use the map
\[ B_*:\Omega^{spin}_*(X) \to H_*(X,\Z): [f:M\to X] \mapsto f_*[M] \]
which is an edge homomorphism in the Atiyah-Hirzebruch spectral sequence.
 Of course, $w=f_*([T^5])=u_*[M]$ is the image under this transformation
 of the  singular
manifold we consider.

If $M$ would admit a pscm, then 
\[ w\in H^+_5(B\pi) \]
Iterated application of theorem \ref{Hom_cond} implies that
\[ 0\ne z\in H^+_2(B\pi) \]
But there is only one two dimensional oriented manifold with positive
curvature, namely $S^2$. Since $\pi_2(B\pi)=0$ any map
$g:S^2\to B\pi$ is null homotopic. In particular $g_*[S^2]=0\in H_2(B\pi,\Z)$,
and therefore $H_2^+(B\pi,\Z)=0$.

This is the desired contradiction and $M$ does not admit a metric with
positive scalar curvature.

The proof relies on the existence of torsion in $\pi$. Therefore, one may
still conjecture that the unstable Gromov-Lawson-Rosenberg conjecture is
true for torsion free groups.

\section{Acknowledgements}
To work on the counterexample
 was inspired by talks of Stephan Stolz where he expressed his opinion that
the original 
GLR-conjecture is false. The author wants to thank Stephan Stolz for
 useful and enlightening
 conversations on the subject. Stolz conjectures that
 a weaker form, the so called stable
GLR-conjecture, is true (compare \cite{Stolz(1994b)}) 
and shows \cite{Stolz(1996b)}  that this
conjecture follows from the Baum-Connes conjecture.

\end{document}